\documentclass[a4paper,10pt]{article}
\usepackage{amsmath,amscd,amsthm,amsfonts,amssymb,epsf,latexsym}
\usepackage[all]{xy}

\makeatletter

\long\def\@makefnt#1{\parindent 1em\noindent
             \hb@xt@1.8em{\hss\@textsuperscript{}}#1}
\long\def\@ftntext#1{\insert\footins{%
     \reset@font\footnotesize
     \interlinepenalty\interfootnotelinepenalty
     \splittopskip\footnotesep
     \splitmaxdepth \dp\strutbox \floatingpenalty \@MM
     \hsize\columnwidth \@parboxrestore
     \color@begingroup
       \@makefnt{%
         \rule\z@\footnotesep\ignorespaces#1\@finalstrut\strutbox}%
     \color@endgroup}}%
\def\subjclass#1{%
   \@ftntext{2010 {\itshape Mathematics Subject Classification. {\rm Primary}
   }\enspace
#1.}}
\def\keywords#1{%
   \@ftntext{{\itshape Keywords.}\enspace #1.}}
\makeatother

\def\moins{\raise 1pt\hbox{{$\scriptstyle -$}}}
\def\plus{\raise 1pt\hbox{{$\scriptstyle +$}} }

\newtheorem{theorem}{Theorem}
\newtheorem{proposition}[theorem]{Proposition}
\newtheorem{lemma}[theorem]{Lemma}
\newtheorem{corollary}[theorem]{Corollary}

\newtheorem{example}[theorem]{Example}

\def\proof{\noindent{\bf Proof.\ }}

\def\qed{~\hbox{$\Box$}}

\def\rank{\mathop{\rm rank}}

\def\Z{{\bf Z}}

\def\Q{{\bf Q}}

\def\\Q{\widetilde{Q}}

\def\qed{~\hbox{$\Box$}}

\begin{document}

\title{\bf A Gysin formula for Hall-Littlewood polynomials}

\author{
Piotr Pragacz\thanks{This work was supported by National Science Center (NCN) grant No. 2014/13/B/ST1/00133.}\\
\small Institute of Mathematics, Polish Academy of Sciences\\
\small \'Sniadeckich 8, 00-656 Warszawa, Poland\\
\small P.Pragacz@impan.pl}



\date{}

\maketitle

\centerline{\it To Bill Fulton on his 75th birthday}


\begin{abstract}
We give a formula for pushing forward the classes of Hall-Littlewood polynomials in Grassmann bundles,
generalizing Gysin formulas for Schur $S$- and $P$-functions.
\end{abstract}

Let $E\to X$ be a vector bundle of rank $n$ over a nonsingular variety $X$
over an algebraically closed field. Denote by $\pi: G^q(E)\to X$ the Grassmann bundle parametrizing rank $q$ quotients of $E$.
Let $\pi_*: A(G^q(E))\to A(X)$ be the homomorphism of the Chow groups of algebraic cycles modulo rational equivalence,
induced by pushing-forward cycles (see \cite[Chap. 1]{Fu1}). There exists an analogous map of cohomology groups.
A goal of this note is to give a formula (see Theorem \ref{main}) for the image via $\pi_*$ of Hall-Littlewood
classes from the Grassmann bundle.

Hall-Littlewood polynomials appeared implicitly in Hall's study \cite{Ha} of the combinatorial
lattice structure of finite abelian $p$-groups, and explicitly in the work of Littlewood on some
problems of representation theory \cite{Li}. A detailed account of the theory of Hall-Littlewood
functions is given in \cite{Mcd}.

The formula in Theorem \ref{main} generalizes some Gysin formulas for Schur $S$- and $P$-functions. 
In particular, it generalizes the formula in \cite[Prop. 1.3(ii)]{Pr5}, and provides
an explanation of its intriguing coefficient. We refer to \cite{FuPr} for general information about the appearance
of Schur $S$- and $Q$-functions in cohomological studies of algebraic varieties.

Let $t$ be an indeterminate. The main formula will be located in $A(X)[t]$, or in the extension $H^*(X,\Z)[t]$ 
of the cohomology ring for a complex variety $X$.
Let $\tau_E: Fl(E)\to X$ be the flag bundle parametrizing flags of quotients of $E$ of ranks $n,n-1,\ldots,1$. 
Suppose that $x_1,\ldots,x_n$ is a sequence of the Chern roots of $E$.
For a sequence $\lambda=(\lambda_1,\ldots,\lambda_n)$ of nonnegative integers, we define
\begin{equation}\label{rle}
R_{\lambda}(E;t)=(\tau_E)_*\bigl(x_1^{\lambda_1}\cdots x_n^{\lambda_n} \prod_{i<j}(x_i-tx_j)\bigr)\,,
\end{equation}
where $(\tau_E)_*$ acts on each coefficient of the polynomial in $t$ separately. (The same convention will be used for
other flag bundles.)

The Grassmann bundle $\pi: G^q(E)\to X$ is endowed with the tautological exact sequence of vector bundles
$$
0\longrightarrow S\longrightarrow \pi^*E\longrightarrow Q \longrightarrow 0\,,
$$
where $\rank(Q)=q$. Let $r=n-q$ be the rank of $S$. Suppose that $x_1\ldots,x_q$ are the Chern roots of $Q$ and 
$x_{q+1},\ldots,x_n$ are the ones of $S$.

\begin{proposition}\label{rl} For sequences $\lambda=(\lambda_1,\ldots,\lambda_q)$ and $\mu=(\mu_1,\ldots,\mu_r)$ of nonnegative integers, we have
$$
\pi_*\bigl(R_{\lambda}(Q;t){R_\mu}(S;t) \prod_{i\le q < j}(x_i-tx_j)\bigr)=R_{\lambda\mu}(E;t),
$$
where $\lambda\mu = (\lambda_1,\ldots,\lambda_q,\mu_1,\ldots,\mu_r)$ is the juxtaposition of $\lambda$ and $\mu$.
\end{proposition}
\proof Consider a commutative diagram
$$
\xymatrix{Fl(Q)\times_{G^q(E)} Fl(S)\ar[d]_{\tau_Q \times \tau_S} \ar[r]^{ \ \ \ \ \ \ \ \ \cong} &Fl(E)\ar[d]^{\tau=\tau_E}\\
          G^q(E)  \ar[r]_{\pi}             &X  }
$$
It follows that
\begin{equation}\label{com}
\pi_*(\tau_Q \times \tau_S)_*=\tau_*\,.
\end{equation}

\medskip

Using Eq.(\ref{rle}) for $Q$ and $S$ and Eq.(\ref{com}), we obtain
$$
\aligned
&\pi_*\bigl(R_{\lambda}(Q;t){R_\mu}(S;t) \prod_{i\le q < j}(x_i-tx_j)\bigr)\\
&=\pi_*\bigl((\tau_Q)_*\bigl(x_1^{\lambda_1}\cdots x_q^{\lambda_q} \prod_{i<j \le q}(x_i-tx_j)\bigr) \cdot (\tau_S)_*\bigl(x_{q+1}^{\mu_1}\cdots x_n^{\mu_r} \prod_{q<i<j}(x_i-tx_j)\bigr) \prod_{i\le q < j}(x_i-tx_j)\bigr)\\
&=\pi_*(\tau_Q \times \tau_S)_*\Bigl(x_1^{\lambda_1}\cdots x_q^{\lambda_q} \prod_{i<j \le q}(x_i-tx_j) x_{q+1}^{\mu_1}\cdots x_n^{\mu_r} \prod_{q < i<j}(x_i-tx_j)\prod_{i\le q < j}(x_i-tx_j)\Bigr)\\
&=\tau_*(x_1^{\lambda_1}\cdots x_q^{\lambda_q}x_{q+1}^{\mu_1} \cdots x_n^{\mu_r} \prod_{i<j}(x_i-tx_j))\\
&=R_{{\lambda}{\mu}}(E;t)\,.
\endaligned
$$
In the argument above, we have used the following equality:
$$
\prod_{i<j \le q}(x_i-tx_j)\prod_{q<i<j}(x_i-tx_j)\prod_{i\le q<j}(x_i-tx_j)=\prod_{i<j}(x_i-tx_j)\,. \qed
$$

We now set
\begin{equation}\label{vm}
v_m(t)=\prod_{i=1}^m \frac{1-t^i}{1-t}=(1+t)(1+t+t^2)\cdots (1+t+\cdots+t^{m-1}).
\end{equation}

Let $\lambda=(\lambda_1,\ldots,\lambda_n)$ be a sequence of nonnegative integers. Consider the maximal subsets  $I_1$,...,$I_d$ in $\{1,\ldots, n\}$, 
where the sequence $\lambda$ is constant.
Let $m_1,\ldots, m_d$ be the cardinalities of $I_1$,...,$I_d$.
So we have $m_1+ \cdots +m_d=n$. We set
\begin{equation}\label{vl}
v_{\lambda}(t)=\prod_{i=1}^d v_{m_i}(t)\,.
\end{equation}

Let $S_n$ be the symmetric group of permutations of $\{1,\ldots,n\}$. We define a subgroup $S_n^{\lambda}$ of $S_n$ as the stabilizer
of $\lambda$. Of course,
$$
S_n^{\lambda}= \prod_{i=1}^d S_{m_i}\,.
$$
Finally, we associate to a sequence $\lambda$ a $(d-1)$-step flag bundle (with steps of lengths $m_i$) 
$$
\eta_{\lambda}: Fl_{\lambda}(E)\to X\,,
$$ 
parametrizing flags of quotients of $E$ of ranks 
\begin{equation}\label{ranks}
n-m_d,n-m_d-m_{d-1},\ldots,n-m_d-m_{d-1}-\cdots-m_2\,.
\end{equation}

\begin{example}\label{nuo} \rm
Let $\nu=(\nu_1>\ldots>\nu_k>0)$ be a strict partition (see \cite[I,1,Ex.9]{Mcd}) with $k\le n$. Let $\lambda=\nu 0^{n-k}$ be the sequence $\nu$ with $n-k$ zeros
added at the end. Then $d=k+1$, $(m_1,\ldots,m_d)=(1^k, n-k)$, $v_{\lambda}(t)=v_{n-k}(t)$, $S^{\lambda}_n=(S_1)^k \times S_{n-k}$, 
and $\eta_{\lambda}: Fl_{\lambda}(E)\to X$ is the flag bundle, often denoted by $\tau_E^k$, parametrizing quotients of $E$ of ranks $k,k-1,\ldots,1$.

If $\lambda=(a^pb^{n-p})$, then $d=2$, $(m_1,m_2)=(p,n-p)$, $v_{\lambda}(t)=v_p(t)v_{n-p}(t)$, $S^{\lambda}_n=S_p\times S_{n-p}$, and $\eta_{\lambda}$ is here the Grassmann bundle $\pi:G^p(E)\to X$.

\end{example}

\smallskip

We shall now need some results from \cite[III]{Mcd}. Let $y_1,\ldots,y_n$ and $t$ be independent
indeterminates. We record the following equation from \cite[III, (1.4)]{Mcd}: 

\begin{lemma}\label{sum} We have
$$
\sum_{w\in S_n} w\Bigl(\prod_{i<j}\frac{y_i-ty_j}{y_i-y_j}\Bigr)=v_n(t)\,.
$$
\end{lemma}

For a sequence $\lambda=(\lambda_1, \ldots, \lambda_n)$ of nonnegative integers, we define
$$
R_{\lambda}(y_1,\ldots,y_n;t)=\sum_{w\in S_n} w\Bigl(y_1^{\lambda_1}\cdots y^{\lambda_n}\prod_{i<j} \frac{y_i-ty_j} {y_i-y_j}\Bigr)
$$ 

Arguing as in \cite[III (1.5)]{Mcd}, we show with the help of Lemma \ref{sum} the following result.
\begin{proposition}\label{div} The polynomial $v_{\lambda}(t)$ divides $R_{\lambda}(y_1,\ldots,y_n;t)$, and we have
$$
R_{\lambda}(y_1,\ldots, y_n;t)=v_{\lambda}(t) \sum_{w\in {S_n/S_n^{\lambda}}} w \Bigl(y_1^{\lambda_1}\cdots y_n^{\lambda_n} \prod_{i<j, 
\lambda_i\ne \lambda_j} \frac{y_i-ty_j}{y_i-y_j}\Bigr)\,.
$$
\end{proposition}

Let us invoke the following description of the Gysin map for the flag bundle $\eta_{\lambda}: Fl_{\lambda}(E)\to X$
with the help of a symmetrizing operator. Recall that $A(Fl_{\lambda}(E))$ as an $A(X)$-module is generated by 
$S_n^{\lambda}$-invariant polynomials in the Chern roots of $E$ (see \cite[Thm 5.5]{BGG}). 
We define for an $S_n^{\lambda}$-invariant polynomial $f=f(y_1,\ldots,y_n)$,
$$
\partial_{\lambda}(f)=\sum_{w\in S_n/S_n^{\lambda}} w\Bigl(\frac{f(y_1,\ldots,y_n)}{\prod_{i<j, \lambda_i\ne \lambda_j} (y_i-y_j)}\Big)\,.
$$
The following result is a particular case of \cite[Prop. 2.1]{Br} (in the situation of Corollary \ref{pno}, the result was shown already in 
\cite[Sect. 2]{Pr}).
 \begin{proposition}\label{gmgr} With the above notation, we have
$$
({\eta_{\lambda}})_*\bigl(f(x_1,\ldots,x_n)\bigr)= \bigl((\partial_{\lambda} f)(y_1,\ldots,y_n)\bigr) (x_1,\ldots, x_n)\,.
$$
\end{proposition}

It follows from Propositions \ref{div} and \ref{gmgr} that
$$
R_{\lambda}(E;t)=v_{\lambda}(t)(\eta_{\lambda})_*\Bigl(x_1^{\lambda_1}\cdots x_n^{\lambda_n} \prod_{i<j, \lambda_i\ne \lambda_j}(x_i-tx_j)\Bigr)\,,
$$
where $x_1,\ldots ,x_n$ are the Chern roots of $E$. 

\smallskip

Let $\lambda$ be a sequence of nonnegative integers. Extending \cite[III, 2]{Mcd}, we set 
\begin{equation}\label{ple}
P_{\lambda}(E;t)= \frac{1} {v_{\lambda}(t)} R_{\lambda}(E;t)\,.
\end{equation}
It follows from Proposition \ref{div} that $P_{\lambda}(E;t)$ is a polynomial in the Chern classes of $E$ and $t$.

Let us record the following particular case.
\begin{corollary}\label{pno} \rm
Let $\nu$ be a strict partition with length $k\le n$. Set $\lambda=\nu 0^{n-k}$. We have
$$
P_{\lambda}(E;t)=(\tau_E^k)_*\Bigl(x_1^{\nu_1}\cdots x_k^{\nu_k} \prod_{i<j,i\le k}(x_i-tx_j)\Bigr)\,.
$$
\end{corollary}

As a consequence of Propositions \ref{rl} and \ref{div}, using Eq.(\ref{ple}), we obtain the following result.
\begin{theorem}\label{main} Let $\lambda=(\lambda_1,\ldots,\lambda_q)$ and $\mu=(\mu_1,\ldots,\mu_r)$ be sequences of nonnegative integers. Then we have
$$
\pi_*\Bigl(\prod_{i \le q < j}(x_i-tx_j) P_{\lambda}(Q;t) P_{\mu}(S;t)\Bigr)=\frac {v_{{\lambda}\mu}(t)} {v_{\lambda}(t)v_{\mu}(t)}  P_{{\lambda}{\mu}}(E;t)\,.
$$
\end{theorem}

We first consider the specialization $t=0$.
\begin{example} \rm 
We recall Schur $S$-functions. 
Let $s_i(E)$ denotes the $i$th complete symmetric function in the roots $x_1,\ldots,x_n$, given by
$$
\sum_{i\ge 0} s_i(E)=\prod_{j=1}^n\frac {1}{1-x_j}\,.
$$
Given a partition 
$
\lambda=(\lambda_1 \ge \ldots \ge \lambda_n\ge 0)\,,
$
we set
$$
s_\lambda(E)= \bigl| s_{\lambda_{i}-i+j}(E) \bigr|_{1\leq i,j \le n}\,.
$$
(See also \cite[I, 3]{Mcd}.)
Translating the Jacobi-Trudi formula ({\it loc.cit.}) to the Gysin map for $\tau_E: Fl(E)\to X$ (see, e.g. \cite[Sect. 4]{Pr5}),
we have
$$
s_{\lambda}(E)=(\tau_E)_*(x_1^{\lambda_1+n-1}\cdots x_n^{\lambda_n}).
$$
We see that $P_{\lambda}(E;t)=s_{\lambda}(E)$ for $t=0$. Under this specialization, the theorem becomes
$$
\aligned
\pi_*\bigl((x_1\cdots x_q)^r s_{\lambda}(Q) s_\mu(S)\bigr)&= \pi_*\bigl(s_{\lambda_1+r,\ldots,\lambda_q+r}(Q) s_{\mu}(S)\bigr)\\
&=s_{\lambda \mu}(E)\,,
\endaligned
$$
a result obtained originally in \cite[Prop. p. 196]{La} and \cite[Prop. 1]{JLP}. 

\smallskip

If a sequence $\lambda=(\lambda_1,\ldots,\lambda_n)$ is not a partition, then $s_{\lambda}(E)$
is either $0$ or $\pm s_{\mu}(E)$ for some partition $\mu$. 
One can rearrange
$\lambda$ by a sequence of operations $(\ldots,i,j,\ldots) \mapsto
(\ldots,j-1,i+1,\ldots)$ applied to pairs of successive integers. Either
one arrives at a sequence of the form $(\ldots,i,i+1,\ldots)$, in which
case $s_{\lambda}(E)=0$, or one arrives in $d$ steps at a partition
$\mu$, and then $s_{\lambda}(E)=(-1)^d s_{\mu}(E)$.

\end{example} 

\begin{corollary}\label{ns} \rm 
Let $\nu$ and $\sigma$ be strict partitions of lengths $k\le q$ and $h\le r$. 
It follows from Eq.(\ref{vm}) that 
$$
\frac {v_{\nu 0^{q-k}\sigma 0^{r-h}}(t)} {v_{\nu 0^{q-k}}(t)v_{\sigma 0^{r-h}}(t)}=\left[{n-k-h} \atop {q-h}\right](t) \cdot (1+t)^e\,,
$$
the Gaussian polynomial times $(1+t)^e$ where $e$ is the number of common parts of $\nu$ and $\sigma$.
Thus the theorem applied to the sequences $\lambda=\nu 0^{q-k}$ and $\mu =\sigma 0^{r-h}$ yields the following equation:
\begin{equation}\label{nu}
\pi_*\bigl(\prod_{i \le q < j}(x_i-tx_j) P_{\nu}(Q;t) P_{\sigma}(S;t)\bigr)= \left[{n-k-h} \atop {q-h}\right](t) \cdot (1+t)^e \cdot P_{{\lambda}{\mu}}(E;t)\,.
\end{equation}
\end{corollary}

We need the following property of Gaussian polynomials, which should be known but we know no precise reference.

\begin{lemma}\label{gauss} At $t=-1$, the Gaussian polynomial 
$$
\left[{a+b} \atop {a}\right](t)
$$ 
specializes to zero if $a b$ is odd and to the binomial coefficient 
$$
\left(\lfloor (a+b)/2\rfloor \atop \lfloor a/2\rfloor \right)
$$ 
otherwise.
\end{lemma}
\proof
We have
$$
\left[{a+b} \atop {a}\right](t)=\frac{(1-t)(1-t^2)\cdots (1-t^{a+b})}{(1-t)\cdots(1-t^a)(1-t)\cdots (1-t^b)}\,.
$$
Since $t=-1$ is a zero with multiplicity 1 of the factor $(1-t^d)$ for even $d$, and a zero with multiplicity 0 for odd $d$,
the order of the rational function $\left[{a+b} \atop {a}\right](t)$ at $t=-1$ is equal to 
\begin{equation}\label{ab}
\lfloor(a+b)/2\rfloor-\lfloor a/2\rfloor-\lfloor b/2\rfloor\,.
\end{equation}
The order (\ref{ab}) is equal to 1 when $a$ and $b$ are odd, and 0 otherwise. In the former case, we get the claimed vanishing, and 
in the latter one, the product of the factors with even exponents is equal to
$$
\left[\lfloor{a+b}/2\rfloor \atop \lfloor a/2 \rfloor \right](t^2)\,.
$$
The value of this function at $t=-1$ is equal to
$
\left[\lfloor {a+b}/2\rfloor \atop \lfloor a/2\rfloor \right](1)
$ 
which is the binomial coefficient 
$$
{\left(\lfloor (a+b)/2\rfloor \atop \lfloor a/2\rfloor \right)}\,.
$$
This is the requested value since the remaining factors with an odd exponent give 2 in the numerator and the same number 
in the denominator. 

The assertions of the lemma follow. \qed 

\medskip

We now consider the specialization $t=-1$. 
\begin{example} \rm
Consider Schur $P$-functions $P_{\lambda}(E)=P_{\lambda}$ (or $P_{\lambda}(y_1,\ldots,y_n)=P_{\lambda}$) defined as follows. For a strict partition 
$\lambda=(\lambda_1 >\ldots >\lambda_k>0)$ with odd $k$,
$$
P_{\lambda}=P_{\lambda_1} P_{\lambda_2,\ldots,\lambda_k}-P_{\lambda_2} P_{\lambda_1,\lambda_3,\ldots,\lambda_k}+\cdots +P_{\lambda_k} P_{\lambda_1,\dots,\lambda_{k-1}}\,,
$$
and with even $k$,
$$
P_{\lambda}=P_{\lambda_1,\lambda_2} P_{\lambda_3,\ldots,\lambda_k}-P_{\lambda_1,\lambda_3} P_{\lambda_2,\lambda_4,\ldots,\lambda_k}+\cdots +P_{\lambda_1,\lambda_k} P_{\lambda_2,\dots,\lambda_{k-1}}\,.
$$
Here, $P_i=\sum s_\mu$, the sum over all hook partitions $\mu$ of $i$, and for positive $i>j$ we set
$$
P_{i,j}=P_i P_j+2\sum_{d=1}^{j-1} (-1)^d P_{i+d}P_{j-d}+(-1)^jP_{i+j}\,. 
$$
(See also \cite[III, 8]{Mcd}.) It was shown in \cite[p. 225]{Sch} that for a strict partition $\lambda$ of length $k$,
$$
P_{\lambda}(y_1,\ldots,y_n)=\sum_{w\in S_n/(S_1)^k\times S_{n-k}} w \Bigl(y_1^{\lambda_1}\cdots y_n^{\lambda_n} \prod_{i<j,i\le k} \frac{y_i+y_j} {y_i-y_j}\Bigr)
$$
(see also \cite[III, 8]{Mcd}).
This implies
$$
P_{\lambda}(E)=(\tau_E^k)_*\Bigl(x_1^{\lambda_1}\cdots x_k^{\lambda_k} \prod_{i<j,i\le k}(x_i+x_j)\Bigr)\,.
$$
By Corollary \ref{pno}, we see that $P_{\lambda}(E)=P_{\lambda}(E;t)$ for $t=-1$. 

We now use the notation from Corollary \ref{ns}.
Specializing $t=-1$ in Eq.(\ref{nu}), we get by Lemma \ref{gauss}
$$
\pi_*\bigl(c_{qr}(Q\otimes S)  P_\nu (Q)  P_\sigma (S)\bigr) =d_{\nu,\sigma}  P_{\nu \sigma}(E)\,,
$$
where $d_{\nu,\sigma}=0$ if $(q-k)(r-h)$ is odd and 
$$
d_{\nu,\sigma}= (-1)^{(q-k)h} \left(\lfloor (n-k-h)/2 \rfloor \atop \lfloor(q-k)/2\rfloor \right)
$$ 
otherwise. This result was obtained originally in \cite[Prop. 1.3(ii)]{Pr5} in a different way.
The present approach gives an explanation of the intriguing coefficient $d_{\nu,\sigma}$.

Suppose that $\lambda=(\lambda_1,\ldots,\lambda_k)$ is not a strict partition.
If there are repetitions of elements in $\lambda$, then $P_{\lambda}$ is zero; 
if not then
$
P_{\lambda}=(-1)^l P_{\mu}\,,
$
where $l$ is the length of the permutation which rearranges
$(\lambda_1,\ldots,\lambda_k)$ into the corresponding
strict partition $\mu$.
\end{example}

\smallskip

We thank Witold Kra\'skiewicz, Itaru Terada and Anders Thorup for helpful discussions, and the referee 
for suggesting several improvements of the text.

\end{document}